\documentclass[10pt]{amsart}

\usepackage{amssymb,amsmath,latexsym}
\usepackage{charter,eucal}
\usepackage{amssymb}
\usepackage{amscd}
\usepackage{hyperref,mathrsfs}

\newtheorem{theorem}{Theorem }[section]

\newtheorem{lemma}[theorem]{Lemma}

\newtheorem{remark}[theorem]{Remark}
\newtheorem{corollary}[theorem]{Corollary}
\newtheorem{proposition}[theorem]{Proposition}
\newtheorem{observation}[theorem]{Observation}
\newtheorem{conjecture}[theorem]{Conjecture}

\newcommand{\bt}{\begin{theorem}}
\newcommand{\et}{\end{theorem}}
\newcommand{\bc}{\begin{corollary}}
\newcommand{\bl}{\begin{lemma}}
\newcommand{\ec}{\end{corollary}}
\newcommand{\el}{\end{lemma}}
\newcommand{\bo}{\begin{observation}}
\newcommand{\eo}{\end{observation}}
\newcommand{\bp}{\begin{proposition}}
\newcommand{\ep}{\end{proposition}}
\newcommand{\br}{\begin{remark}}
\newcommand{\er}{\end{remark}}
\newcommand{\brt}{\begin{result}}
\newcommand{\ert}{\end{result}}
\newcommand{\eop}{\hspace*{\fill}{\footnotesize $\blacksquare$}}

\def\I{\mathbf{I}}
\def\nI{\not\mathbf{I}}

\def\PG{\mathbf{PG}}

\newcommand{\Aut}{\mathrm{Aut}}

\newcommand{\hP}{\mathbf{P}}
\newcommand{\K}{\mathbb{K}}
\newcommand{\proj}{\mathrm{proj}}

\newcommand{\mU}{\mathcal{U}}
\newcommand{\mF}{\mathbf{F}}
\newcommand{\hF}{\mathbf{F}}

\newcommand{\he}{\mathrm{id}}
\newcommand{\mP}{\mathcal{P}}
\newcommand{\mB}{\mathcal{B}}
\newcommand{\mL}{\mathcal{L}}
\newcommand{\mQ}{\mathcal{Q}}
\newcommand{\mN}{\mathcal{N}}
\newcommand{\mO}{\mathcal{O}}
\newcommand{\mA}{\mathcal{A}}
\newcommand{\mG}{\mathcal{G}}

\newcommand{\mS}{\mathcal{S}}
\newcommand{\mE}{\mathcal{E}}
\newcommand{\F}{\mathcal{F}}

\newcommand{\bA}{\mathbb{A}}
\newcommand{\bP}{\mathbb{P}}
\newcommand{\bT}{\mathbb{T}}
\newcommand{\vo}{\overline{0}}

\title{An obstruction relating locally finite polygons to translation quadrangles}


\author{Koen Thas}

\thanks{}

\address{Ghent University, Department of Mathematics,\\
Krijgslaan 281, S25, B-9000 Ghent, Belgium}

\email{koen.thas@gmail.com}

\date{}

\begin{document}

\maketitle

\begin{abstract}
One of the most fundamental open problems  in Incidence Geometry, posed by Tits in the 1960s, asks for the existence of so-called ``locally finite
generalized polygons'' | that is, generalized polygons  with ``mixed parameters'' (one being finite and the other not).
In a more specialized context, another long-standing problem (from the 1990s) is as to whether
the endomorphism ring of any translation generalized quadrangle is a skew field (the answer of which is known in the finite case). 
(The analogous problem for projective planes, and its positive solution, the ``Bruck-Bose construction,'' lies at the very base of the whole theory of translation planes.)\\
In this short note, we introduce a category, representing certain very specific embeddings of generalized polygons, which 
surprisingly controls the solution of both (apparently entirely unrelated) problems.   
\\ 
\end{abstract}

\bigskip
\tableofcontents

\bigskip
\section{Introduction}

\subsection{First problem | locally finite polygons}
Consider a generalized $n$-gon $\Gamma$ with $s+1$ points on each line and $t+1$ lines through each point and let $s,t>1$ (and note that $st$ is allowed to be infinite) | its incidence graph is a bipartite graph of diameter $n$ and girth $2n$. If $n$ is odd, then it is easy to show that $s=t$, see \cite[1.5.3]{POL}. If $n$ is even, though, there are examples where $s\neq t$, a most striking example being $n=8$ in which case a theorem of Feit and Higman \cite{FeHi} implies that if $st$ is finite, $2st$ is a perfect square and so $s$ is never equal to $t$. If both $s$ and $t$ are finite, they are bounded by each other; to be more specific, $s\leq t^2\leq s^4$ for $n=4$ and $n=8$ by results of Higman \cite{Hi75} (1975),
and $s\leq t^3\leq s^9$ for $n=6$ by Haemers and Roos \cite{HR} (1981). These results can also be found in \cite{POL}, \S 1.7.2.
For other even values of $n$, $\Gamma$ cannot exist by a famous result of Feit and Higman which also appeared  in their 1964 paper \cite{FeHi}. In {\em loc. cit.}, 
several necessary divisibility conditions involving the parameters $s$ and $t$ of a generalized $n$-gon can be found (with $n \in \{ 3,4,6,8\}$), given the existence of such an object with these parameters. 

An old and notorious question, first posed by Jacques Tits in the 1960s, represents the largest gap in our knowledge about parameters of generalized polygons:

\begin{quote}
``Do there exist {\em locally finite} generalized polygons?''
\end{quote}

In other words, do there exist, up to duality, (thick) generalized polygons with a finite number of points incident with a line,
and an infinite number of lines through a point?\footnote{The aforementioned question can be found as Problem 5 in the ``Ten Most Famous Open Problems'' chapter of 
Van Maldeghem's book \cite{POL}, see also \S 10 of \cite{JATGP}, etc.}
(Note that in Van Maldeghem's book \cite{POL}, such generalized polygons are called {\em semi-finite}.)

There is only a very short list of results on Tits's question. All of them comprise the case $n = 4$.
\begin{itemize}
\item
P. J. Cameron \cite{PJC} showed in 1981 that if $n=4$ and $s=2$, then $t$ is finite. 
\item
In \cite{AEB} A. E. Brouwer shows the same thing for $n=4$ and $s=3$ and the proof is purely combinatorial (unlike a nonpublished but earlier proof of Kantor \cite{Kanunp}). 
\item
More recently, G. Cherlin used Model Theory (in \cite{Cher}) to handle the generalized $4$-gons  with five points on a line.
\item
For other values of $n$ and $s$ (where $n$ of course is even), nothing is known without any extra assumptions. 
\end{itemize}

Apart from the aforementioned results, there is only one other ``general result'' on parameters of generalized polygons (so without invoking additional structure through, e.g., the existence of certain substructures or the occurrence of certain group actions):

\bt[Bruck and Ryser \cite{BR}, 1949]
If $\Gamma$ is a finite projective plane of order $m$, $m \in \mathbb{N}^\times$, and $m \equiv 1,2\mod{4}$, then $m$ is the sum of two perfect squares.
\et

\subsection{Second problem | endomorphisms of translation polygons}
In Incidence Geometry, one is often concerned with categories $\mathcal{C}$ of abelian groups such that if $A$ is an object of $\mathcal{C}$,
and $\mathrm{Hom}(A,A)$ is the endomorphism ring of $A$, where morphisms in $\mathcal{C}$ are group homomorphisms that
preserve some specified substructure of the objects,
then $\mathrm{Hom}(A,A)$ is a skew field, and $A$ is a left or right vector space over $\mathrm{Hom}(A,A)$. The reason of this interest is that 
such results lead to projective representations of associated geometrical objects (cf. the next paragraph).

One such classical example is the category $\mathcal{TP}$ of translation groups of affine (or projective) translation planes, where morphisms are
prescribed to preserve the plane attached to it, and all ``points at infinity.'' An early result in this theory is that, indeed, all $\mathrm{Hom}(A,A)$s (which are called {\em kernels}) are 
skew fields, and this leads to the fact that one can represent the associated translation planes in a projective space. This representation is the 
famous  ``Bruck-Bose'' representation, and is arguably (?) the most fundamental tool to study translation planes.\\

A category which is related but much more difficult and mysterious, is that of translation groups of generalized quadrangles, say $\mathcal{TGQ}$.
In the same way as for planes one defines ``kernels,'' and again one aims at proving that all kernels are skew fields.
So the challenge is, very roughly, the question {\em whether all translation quadrangles can be embedded in projective space over a skew field.}

In the finite case, where the answer was proved to be positive already in the 1980s, this embedding result became the most basic tool of the theory, such as for the planes, see \cite{TGQ}.


More precisely,
let $\Gamma^x$ be an infinite translation generalized quadrangle (TGQ) (see the next section for the definitions), with translation group $T$. Let $z$ be an affine point, that is, a point not collinear with $x$, and $\{U\}_z$ be the set of lines incident with $z$.
If $V$ is any such line, put $v := \proj_Vx$. Define $\mathbb{K}$ as the set of endomorphisms of $T$ that map every $T_V$ into itself. Then $\mathbb{K}$ is a ring (with multiplicative identity) without zero divisors \cite{Jos}, and for any such $V$, we have that $T$, $T_V$ and $T_v$ are left $\mathbb{K}$-modules.

\begin{conjecture}[Linearity for TGQs]
$\mathbb{K}$ is a skew field.
\end{conjecture}

The property that $\mathbb{K}$ is a skew field allows one to interpret $T$, $T_V$ and $T_v$ ($V \in \{ U\}_z$) as vector spaces
over $\mathbb{K}$, so as to represent $\Gamma^x$ in projective space (over $\mathbb{K}$), just as in the finite case
| see \cite[Chapter 8]{PT}.  Throughout this paper, we will call TGQs for which the kernel is a skew field {\em linear}.

\br{\rm
Note that both projective planes and generalized quadrangles are special cases of the class of generalized $n$-gons (respectively the case $n = 3$ and $n = 4$). ``Translation generalized $n$-gons'' only make sense for these two cases, though (cf. \cite{POL}).
}\er

\br{\rm
For some special cases the aforementioned conjecture has been proved to be true; planar TGQs \cite{Jos2}, TGQs with a strongly regular translation center \cite{Jos2} and of course finite TGQs \cite{PT} all satisfy the conjecture.
Also, it has been shown in \cite{Jos3} that the more restricted ``topological kernel'' of a compact connected topological TGQ is a skew field.\footnote{Due to these rather restricted partial results, it is usually {\em assumed} that the kernel of a TGQ be a skew field.}}
\er

\br
{\rm
In \cite{Jos} it is claimed (in Corollary 3.11) that the kernel of a TGQ always is a skew field. In the proof however, the author uses his Proposition 3.10 which states 
that any three distinct lines on the translation point, together with any affine point, generate a plane-like
subGQ. This result is not true (even not in the finite case) | in fact, only a very restricted class of TGQs
has this property. (Still, Theorem 3.11 of {\em loc. cit.} shows that if a TGQ {\em does} satisfy this property,
it indeed is linear.) The paper \cite{Jos} contains many other interesting results on infinite TGQs.
}
\er

\subsection{The present paper}

In this paper, it is our intention to show that  if {\em one} curious property could be understood,   
these two seemingly (very) unrelated problems would follow. (In both problems, one wants to control the parameters of the 
polygon, but this is barely evidence for such connections to exist; see also \cite{Order} for an elaborate discussion.)

We will show that if  there do not exist pointed generalized polygons $(\mU,u)$, $(\mU',u)$ (where $u$ is a point), together with an isomorphism
$\gamma: \mU \mapsto \mU'$ such that $\mU'$ is ideally embedded in $\mU$, and $u \in \mU' \leq \mU$ is fixed linewise by $\gamma$, 
then both problems have a positive answer (cf. Theorem \ref{MRE} below). 

As such, this property represents a ``common geometrically divisor'' of both problems.

\section{Some definitions}

\subsection{}
In this paper,
a {\em generalized $n$-gon}, $n \in \mathbb{N} \setminus \{0,1,2\}$, is a point-line incidence geometry $\Gamma = (\mP,\mB,\I)$ for which  the following axioms are satisfied:

\begin{itemize}
\item[(i)]
$\Gamma$ contains no ordinary $k$-gon (as a subgeometry), for $2 \leq k < n$;
\item[(ii)]
any two elements $x,y \in \mP \cup \mB$ are contained in some ordinary $n$-gon in
$\Gamma$;
\item[(iii)]
there exists an ordinary $(n + 1)$-gon in $\Gamma$.\\
\end{itemize}

A {\em generalized polygon} (GP) is a {\em generalized $n$-gon} for some $n$, and $n$ is called the {\em gonality} of the GP. 
Note that projective planes are generalized $3$-gons. Generalized $4$-gons are also called {\em generalized quadrangles}.
By (iii), generalized polygons have at least three points per line and three lines per point. So by this definition, we do not consider ``thin'' polygons.\footnote{For thin polygons, the questions under consideration obviously make no sense.}
Note that points and lines play the same role; this is the principle of ``duality''.

A {\em pointed GP} is a pair $(\Gamma,r)$ where $\Gamma$ is a GP and $r$ is either a point or a line in $\Gamma$.

It can be shown that generalized polygons have an {\em order}; there exist constants $s, t$ such that 
the number of points incident with a line is $s + 1$, and the number of lines incident  with a point is 
$t + 1$, cf. \cite[1.5.3]{POL}. If the gonality is odd, then $s = t$, cf. {\em loc. cit.}

\br{\rm
Generalized polygons were introduced by Tits in a famous work on triality \cite{Ti} of 1959,  in order to propose an axiomatic and combinatorial treatment for  semisimple algebraic groups (including Chevalley groups  and groups of Lie type) of relative rank 2.}
\er

\subsection{}
A {\em sub generalized polygon} or {\em subpolygon} of a GP $\Gamma = (\mP,\mB,\I)$ is a GP $\Gamma' = (\mP',\mB',\I')$ for which $\mP \subseteq \mP$,
$\mB' \subseteq \mB$ and $\I' \subseteq \I$. A subpolygon has the same gonality as its ``ambient polygon'' $\Gamma$. It is {\em full} if for any line $L$ of $\Gamma'$, we have that $x \I'  L$ if and only if $x \I L$. Dually, we define {\em ideal subpolygons}.

\subsection{}
An {\em automorphism} of a generalized polygon $\Gamma = (\mP,\mB,\I)$ is a bijection of $\mP \cup \mB$ which preserves $\mP$, $\mB$ and incidence.
The full set of automorphisms of a GP forms a group in a natural way | the {\em automorphism group} of $\Gamma$, denoted $\Aut(\Gamma)$. It is one of its most important invariants.
If $B$ is an automorphism group of a generalized polygon $\Gamma = (\mP,\mB,\I)$, and $R$ is a subset of $\mP$, $B_{[R]}$ is the subgroup of $B$ fixing $R$ pointwise (in this notation, a line is also considered to be a point set).

{\em Morphisms} between GPs are defined similarly as automorphisms, and if $A$ and $B$ are GPs, $\mathrm{End}(A,B) = \mathrm{Hom}(A,B)$ denotes the set of all morphisms  $A \mapsto B$. If $A = (L,l)$ and $B = (L',l')$ are pointed GPs, elements of $\mathrm{End}((L,l),(L',l'))$ map $l$ to $l'$.
Also, if $(L,l) = (L',l')$, $\mathrm{End}((L,l),(L,l))$ is shortened to $\mathrm{End}(L,l)$, and $\mathrm{End}(L,[l])$ denotes the subset
of elements that fix $l$ {\em elementwise}.

\subsection{}
A {\em translation generalized quadrangle} (TGQ) \cite{PT,TGQ} $\Gamma^x$ is a generalized quadrangle for which there  is an abelian automorphism group $T$ that fixes each line incident with the point $x$, while acting sharply transitively on the points not collinear with $x$ (= the {\em affine points} ``w.r.t. $x$'').

Suppose $z$ is an affine point (w.r.t. $x$).
Let $\F = \{T_M \vert M \I z\}$ and $\F^* = \{T_m \vert x\sim m \sim z\}$; if $m \I M$, we 
denote $T_m$ also by $T_M^*$.
Then for all $L \I z$ we have:

\begin{itemize}
\item
$T_L \leq T_L^* \ne T_L$;
\item
$T_LT_M^* = T$ for $M \ne L$;
\item
$T_AT_B \cap T_C = \{\he\}$ for distinct lines $A,B,C$.
\item
$\{T^*_L/T_L\} \cup \{T_LT_M/T_L \vert  M \ne L\}$ is a partition of $T/T_L$.
\end{itemize}

Call $(T,\F,\F^*)$ a {\em Kantor family} of $\Gamma^x$.

Conversely, from families $\F$ and $\F^*$ with these properties in an abelian group $T$, one can construct 
a TGQ for which $(T,\F,\F^*)$ is a Kantor family, using a natural group coset geometry respresentation \cite{PT,TGQ}. If one starts from a TGQ as above, it is isomorphic to the reconstructed coset geometry (the isomorphism class of the geometry is independent of the chosen affine point $z$). This representation method was noted for a more general class of GQs (namely, for so-called ``elation generalized quadrangles'') by
Kantor \cite{Ka} in the finite case, and carries over without much change to the infinite case \cite{BaPa}.

\section{Statement of main result}
\label{CATE}

Define the category $\mE$ as follows: its objects are pairs $(X,\eta)$, where
$X = (\Gamma,x)$ is a pointed GP, and $\eta \in \mathrm{End}(\Gamma,[x])$ is injective but not bijective.
If $X = ((\Gamma,x),\eta)$ and $Y = ((\Gamma',x'),\eta')$ are objects in $\mE$, $\mathrm{Hom}(X,Y)$
consists of morphisms $\beta: \Gamma \mapsto \Gamma'$ sending $x$ to $x'$ (``pointed morphisms'') such that
$\beta \circ \eta = \eta' \circ \beta$.

Note that it follows that $\Gamma^\eta \cong \Gamma$, that $x \in \Gamma^\eta$, and that $\Gamma^\eta$ is ideally embeded
in $\Gamma$. The isomorphism $\eta$ describes the embedding, and fixes $x$ linewise.

\bt[Main result]
\label{MRE}
If $\mE$ is empty, then locally finite polygons do not exist, and all TGQs are linear.
\et 

In fact, only the ``locally finite part'' of $\mE$ and the part for which $X = (\Gamma,x)$ ``is'' a TGQ matter for the proof.

Note that the morphism $\eta: (\Gamma,x) \mapsto (\Gamma,x)$ generated a sequence
\begin{equation}
\cdots \xrightarrow{\eta} (\Gamma,x) \xrightarrow{\eta} \cdots \xrightarrow{\eta} (\Gamma,x) \xrightarrow{\eta} \cdots
\end{equation}
of which the limits in both directions (that is, $\bigcup_{i \in \mathbb{Z}}\Gamma^{\eta^{i}}$ and $\bigcap_{i \in \mathbb{Z}}\Gamma^{\eta^i}$)
are stable under $\eta$.

\section{Obstruction for TGQs}

\subsection{Setting}

$\Gamma = \Gamma^x = (\Gamma^x,T)$ is a TGQ with translation point $x$ and translation group $T$.  We suppose that the number of points (and so also the number of lines) is not finite. The kernel $\K$ is 
defined as above, and we suppose, by way of contradiction, that $\K$ is not a skew field.

\subsection{The GQs $\Gamma(\alpha,z)$}

Let $\alpha \in \K^\times$, and $z$ an affine point. To avoid trivialities, suppose that $\alpha$ is not a unit.
Let $(T,\F,\F^*)$ be the Kantor family in $T$ defined by $z$. 
It is easy to see that $(T^\alpha,\F^{\alpha},{(\F^*)}^\alpha)$ (obvious notation) defines a Kantor family
in $T^\alpha$. Moreover, the GQ $\Gamma(\alpha,z)$ defined by this Kantor family is {\em thick} (it has at least three points per line), and it is  ideal.
One can derive these properties easily from Lemma \ref{inj} below.

Throughout this section, we will fix $z$, so that we write $\Gamma(\alpha)$ instead of $\Gamma(\alpha,z)$.



\subsection{Injectivity}

The proof of the next lemma is essentially the same as in the finite case (which can be found in \cite[Chapter 8]{PT}). We include its proof for the sake of convenience.

\bl
\label{inj}
Each element of $\mathbb{K}$ is injective.
\el
{\em Proof}.\quad
Suppose that $\beta\in \K$ is such that $\ell_{0}^\beta = \he$ for some $\ell_0\in T_0\setminus\{\he\}$, $T_0 \in \{T_U \vert U \I z\}$; then we must show that $\beta = 0$. (The choice of $T_0$ is arbitrary. If $\beta$
has a fixed point not in  $\bigcup_{V \in \{U\}_z} V$, then it has a fixed point in each $V\setminus\{\he\}$ as well.) Assume the contrary. Choose any element $\ell_{i}\in T_{i}\setminus\{\he\}$, with $T_i \ne T_0$.
Then the point $\ell_{0}\ell_{i}$ is at distance two from $\he$ in the collinearity graph of $\Gamma^x$. Since $\Gamma^x$ is thick there exist elements $\ell_l\in T_l\setminus\{\he\}$ and
$\ell_k\in T_k\setminus\{\he\}$, $\{T_l,T_k\}\cap \{T_0,T_1\} = \emptyset$ and $T_l\neq T_k$, such that 
\begin{equation}
\ell_0\ell_i = \ell_l\ell_k.  
\end{equation}
Letting $\beta$ act yields $\ell_{i}^\beta=\ell_{l}^\beta \ell_{k}^\beta$. First suppose that $\ell_{l}^\beta=\he$;  then $\ell_{i}^\beta=\ell_{k}^\beta$. Since $T_i\cap T_k =\{\he\}$ we obtain that
$\ell_{i}^\beta = \he$. Analogously $\ell_k^\beta = \he$ implies that $\ell_i^\beta = \he$. Next suppose that neither $\ell_{l}^\beta$ nor $\ell_{k}^\beta$ equals $\he$. In this case the line $T_l\ell_{k}^\beta$ of $\Gamma^x$ intersects the line $T_k$ in $\ell_{k}^\beta\neq \he$
and intersects the line $T_i$ in $\ell_{i}^\beta\neq \he$. Hence we have found a triangle in $\Gamma^x$, a contradiction. We conclude that $\ell_{i}^\beta=\he$, and henceforth that $V^\beta = \he$, for all $V \in \{ T_U \vert U \I z\}$. By the
connectedness of $\Gamma^x$ we know that $T = \langle V \vert V \in \{T_U \vert U \I z\}\rangle$, and hence it follows that $T^\beta = \he$, that is, $\beta = 0$.  
\eop

\bc
\label{cor}
For each $\alpha \in \mathbb{K} \setminus \{0\}$, $\Gamma(\alpha) \cong \Gamma$.
\eop 
\ec

\bigskip
\subsection{Proof of Theorem \ref{MRE} for TGQs}
Suppose $\Gamma^x$ is a TGQ with endomorphism ring $\mathbb{K}$.
Suppose $\mathbb{K}$ is not a skew field; then there exists a $\zeta \in \mathbb{K}$ which is not invertible. 
By definition of $\mathbb{K}$ and Corollary \ref{cor}, it follows that $((\Gamma,x),\zeta)$ is an object in $\mE$.\\

\section{Obstruction for locally finite polygons}
\label{LFpart}

Let $\Gamma$ be a generalized polygon. An ordered set $\mL$ of lines is {\em indiscernible} if for any
two increasing sequences $M_1,M_2,\ldots,M_n$ and $M_1',M_2',\ldots,M_n'$ (of the same length $n$) of lines of $\mL$, there 
is an automorphism of $\Gamma$ mapping $M_i$ onto $M_i'$ for each $i$. It is indiscernible {\em over $D$}, if $D$ is a finite set of points and lines fixed by the automorphisms just described.

By combining the Compactness Theorem and Ramsey's Theorem \cite{Hodge} (in a theory which has a model in which a given definable set is infinite), one can prove the following.

\bt[G. Cherlin \cite{Cher}]
\label{Cherlem}
Suppose there is an infinite locally finite generalized $n$-gon with finite lines. Then there is an
infinite locally finite generalized $n$-gon $\Gamma$ containing an indiscernible sequence $\mathcal{L}$ of parallel (= mutually skew) lines, of any 
specified order type. The sequence may be taken to be indiscernible over the set $D$ of all points incident with one fixed line $L$ of $\Gamma$.
\et

Clearly, $D$ may supposed to be fixed pointwise in Theorem \ref{Cherlem}.

\br
\label{pol}
{\rm
\begin{itemize}
\item[{\rm(i)}]
Cherlin states the theorem only for generalized quadrangles, but the same statement also holds for ``general'' generalized polygons.
\item[{\rm(ii)}]
In Theorem \ref{Cherlem}, $\Gamma$ may supposed to be generated by $\mathcal{L} \cup \{L\}$ \cite{Cherpriv} (where $L$ is seen as a point set).
\end{itemize}}
\er

\subsection{Setting}

Let $\mL$ be a fixed indiscernible set of (mutually skew) lines of $\Gamma$ over $L$ | here $\Gamma$ is generated by
$\mL \cup \{L\}$ (recall Remark \ref{pol}(ii)).  
We suppose that $\mL$ is (infinitely) countable.

If $S \subseteq \mL$, by $\Gamma(S)$ we denote the full subpolygon generated by $S \cup \{L\}$ (since $L$
is fixed throughout, we do not specify $L$ in the notation). So $\Gamma = \Gamma(\mL)$.

Call a sequence of lines $N_{i_1},N_{i_2},\ldots,N_{i_l}$ in $\mL$ {\em increasing} if $i_l < i_{l'}$ for 
$l < l'$.

\bl
\label{indisc}
If $S \subset \mL$ is a subset of $\mL$, then $\Gamma(S)$ is disjoint from any line of 
$\mL \setminus S$.
\el
{\em Proof}.\quad
Immediate from indiscernibility.
\eop 

 Let $\mathrm{Sub}(\Gamma)$ be the set of subGPs of $\Gamma$ (including thin subGPs).

\bl
The map  $\Psi: 2^{\mL} \mapsto \mathrm{Sub}(\Gamma): S \mapsto \varphi(S) = \Gamma(S)$ 
is an injection. 
\el
{\em Proof}.\quad
Immediate.
\eop

The following lemma is folklore, and easy to prove. (In its statement, ``countable'' also comprises
the finite case.)

\bl
\label{lemcount}
Let $\Delta$ be a generalized polygon, and $K$ be a subset of points of countable size.
Let $\Delta(K)$ be the subpolygon of $\Delta$ generated by $K$. Then the number of points and 
lines of $\Delta(K)$ is countable.
\eop 
\el

\bc
We have that $t$ is countable, as is the number of points and lines of $\Gamma$. 
\ec

{\em Proof}.\quad
Since $\Gamma$ is generated by $\{L\} \cup \mL$, it is generated by a countable number of points.
By Lemma \ref{lemcount}, the number of points and lines is countable. Clearly $t$ also is.
\eop

\subsection{The topology $(\mathrm{Sym}(X),\tau)$}
\label{toptau}

Let $X$ be the point set of $\Gamma$ and let $\mathrm{Sym}(X)$ be the symmetric group on $X$. We endow $\mathrm{Sym}(X)$ with the topology $\tau$ of {\em pointwise convergence} |
a subset of $\mathrm{Sym}(X)$ is {\em closed} if and only if it is the automorphism group of some first order structure. So $\mathrm{Aut}(\Gamma)$ is closed in $\tau$. 
Note that a subset of $\mathrm{Sym}(X)$ is open if and only if it contains the elementwise stabilizer of some finite subset.
The closed sets $C$ are characterized by the following property:

\begin{quote}
($\gamma$)\quad Let $g \in \mathrm{Sym}(X)$. Then $g \in C$ if and only if for every finite subset $A \subset X$ there is a $g_A \in C$ that agrees with $g$ on $A$.
\end{quote} 

In particular, $(\gamma)$ applies to $C = \Aut(\Gamma)$.

\subsection{Automorphisms preserving $\mL$}
\label{subaut}

We index $\mL$ by $\mathbb{Q}$:
$\mL  = \{ M_i \}_{i \in \mathbb{Q}}$, and let $\mathbb{Q}$ be endowed with the natural (linear, dense) order $\leq$. By assumption, $\mL$ is indiscernible
over $L$ w.r.t. $(\mathbb{Q},\leq)$.
Put $\mathrm{Aut}(\Gamma)_{[L]} =: A$. Let $\alpha$ be any order preserving permutation of the index set $\mathbb{Q}$.
Let $\phi: \mathbb{N} \mapsto \mathbb{Q}$ be a bijection, and define, for each $i \in \mathbb{N}$:
\begin{equation}
 \mL_i = \{M_{\phi(j)} \vert 0 \leq j \leq i \}.     
 \end{equation}

So $\mL = \bigcup_i \mL_i$ and $\mL_i \subseteq \mL_j$ for $i \leq j$.
For each $k \in \mathbb{N}$, let $A_k$ be the set of elements of $A$ which have the same action as $\alpha$ on $\mL_k$; it is not an empty set by indiscernibility and 
finiteness of $\mL_k$.

The next lemma is obvious.

\bl
\label{gengeom}
Let $\Delta$ be a GP, and $Y$ a generating set (of points, say). If $\alpha \in \mathrm{Aut}(\Delta)$, its action on $\Delta$ is completely determined by its action on $Y$.
\eop
\el

(If $\alpha' \in \mathrm{Aut}(\Delta)$ would have the same action on $Y$, $\alpha{\alpha'}^{-1}$ fixes $Y$ pointwise, so also the subGP generated by $Y$. But this is 
$\Delta$.)

Since $\{L\} \cup \mL$ generates $\Gamma$, it is (by the previous lemma) clear that if $\bigcap_lA_l =: \overline{A} \ne \emptyset$, its size 
is precisely $1$ | that is, there is a unique automorphism (denoted $\chi$) of $\Gamma$ fixing $L$ pointwise and inducing $\alpha$ on $\mL$.
This is exactly what (i) of the following lemma says.

\bl
\label{partaut1}
\begin{itemize}
\item[{\rm(i)}]
$\overline{A} \ne \emptyset$, and so $\chi \in \mathrm{Aut}(\Gamma)$ is well-defined. 
\item[{\rm(ii)}]
$\chi$ stabilizes $\mL \cup L$.
\end{itemize}
\el
{\em Proof}.\quad
(i) 
Take any $k \in \mathbb{N}$; then the elements of $A_k$ agree on $\Gamma(\mL_k)$ (that is, if $\gamma, \epsilon$ are in $A_k$, then $\gamma\epsilon^{-1} = \mathbf{1}_{\vert A_k}$). Clearly, the elements of $A_{k'}$ also agree on 
$\Gamma(\mL_{k})$ with the elements of $A_k$ if $k' > k$. Now define $\chi$ (inductively) by passing to the limit $k \mapsto \infty$. Since $\Gamma = \bigcup_j\Gamma(\mL_j)$, $\chi$ is a 
well-defined element of $\mathrm{Sym}(X)$ (where $X$ is as in \S \ref{toptau}). Consider any finite subset $R$ of $X$. Then there is some $\ell \in \mathbb{N}$ for which $R \subseteq \Gamma(\mL_\ell)$.
Any element $\beta \in A_\ell$ coincides with $\chi$ on $R$. So since $\mathrm{Aut}(\Gamma)$ is closed, $\chi  \in A$.

(ii) is immediate.  
\eop

(Note that if, for each $i \in \mathbb{N}$, $\chi_i \in A_i$ is chosen arbitrarily, $(\chi_j)_j$ converges to $\chi$ in $(\mathrm{Sym}(X),\tau)$.)

\bt[Determination of $A_{\mL}$]
\label{partaut2}
\begin{itemize}
\item[{\rm(i)}]
$A_\mL$ contains a subgroup $\mO$ isomorphic to the complete group of order preserving permutations of $(\mathbb{Q},\leq)$ acting naturally on $\mathbb{Q}$.
\item[{\rm (ii)}]
$A_{\mL}$ acts $k$-homogeneously on $\mL$ for any $k \in \mathbb{N}^\times$.
\end{itemize}
\et
{\em Proof}.\quad
By Lemma \ref{gengeom}, the kernel of the action of $A_\mL$ on $\mL$ is trivial. Hence (i) follows from Lemma \ref{partaut1}.\\

(ii) follows from (i).
\eop 

\br{\rm
\begin{itemize}
\item[(a)]
It is not hard to show that $(A_\mL,\mL)$ is in fact {\em isomorphic} to the complete group of order preserving permutations of $(\mathbb{Q},\leq)$ acting naturally on $\mathbb{Q}$ (if $A_\mL$ would properly contain the order preserving permutations of $(\mathbb{Q},\leq)$, one can easily show to arrive at a contradiction).
In that case, we can add to (ii) that the action is not $2$-transitive.
\item[(b)]
By (a), $A_{\mL}$ cannot contain involutions.
\end{itemize}
}
\er

\bc
$\mathrm{Aut}(\Gamma)$ is uncountable. As a direct consequence, $\Gamma$ is infinitely generated.
\eop 
\ec

\bigskip
\subsection{Proof of Theorem \ref{MRE} for locally finite polygons}
Now let $I = (a,b)$ be an open interval in $\mathbb{Q}$ with $a < 0 < b$, and let $\zeta \in \mO$ be a homothecy with center $0$ and 
factor $f < 1$. Then $\Gamma(I)^{\zeta} \leq \Gamma(I) \ne \Gamma(I)^{\zeta}$, and $L$ is fixed pointwise by $\zeta$.
So $((\Gamma(I),L),\zeta)$ is an object in $\mathcal{E}$.

This concludes the proof of Theorem \ref{MRE}.\\

\br{\rm
\begin{itemize}
\item[(a)]
Using the stability property referred to in \S \ref{CATE}, the author has recently proved that for a TGQ of order $(s,t)$
to be nonlinear, $s$ and $t$ must be equal, and the characteristic of its kernel then must be $0$. 
\item[(b)]
Currently we are trying to develop the automorphic theory set up in \S \ref{LFpart}.
\end{itemize}
}
\er


\end{document}